\documentclass{article}
\begin{document}

\begin{center}
\large \textbf{Set-Stabilizers in Solvable Permutation Groups}
\end{center}

\bigskip
\bigskip

\begin{center}
\textbf{ABSTRACT}
\end{center}

\medskip
\noindent Let $G$ be a finite solvable permutation group. Then 
modulo a possibly trivial normal elementary abelian 3-subgroup,
some set-stabilizer in $G$ is a 2-group.

\bigskip
\bigskip
\bigskip
\bigskip
\bigskip

\noindent David Gluck\newline
Department of Mathematics\newline
Wayne State University\newline
Detroit, Michigan 48202\newline
d.gluck@me.com

\bigskip
\bigskip
\bigskip
\bigskip
\bigskip
\bigskip
\bigskip
\bigskip
\bigskip
\bigskip


\bigskip
\bigskip
\bigskip
\bigskip


\newpage

\begin{center}
\textit{ To the memory of Marty Isaacs}
\end{center}
\bigskip

\begin{center}
\textbf{1. Introduction}
\end{center}
\bigskip

Let $(G,\Omega)$ be a finite solvable permutation group.
We want to find a subset of  $\Omega$ whose setwise stabilizer in $G$ 
is ``small''. The best we can hope for is a subset with trivial stabilizer in
$G$, or equivalently, a set which lies in a regular $G$-orbit on the power
set of $\Omega$.

Loosely speaking, the prime 2 is the obstruction to the existence of
regular orbits on the power set. Indeed, if $|G|$ is odd, then $G$ has a 
regular orbit on the power set by [2, Corollary 1], but a Sylow 2-subgroup
of the symmetric group on $2^n$ points has no regular orbit on the power 
set of $\{1,2,\ldots, 2^n\}$. One might guess, therefore, that there always exists
a subset of $\Omega$ whose stabilizer in $G$ is a 2-group. However, this is 
false; in the doubly transitive solvable group of degree 8 and order 168,
every subset of the permuted set is stabilized by an element of order 3.

Our main result is Theorem A, a corollary of Theorem 3.2. Theorem A shows that the guess above is, nevertheless, not 
far from the truth. Theorem A is a considerable improvement over the earlier
result ([3, Lemma 7], [7, Corollary 5.7a]) that some set-stabilizer must be a $\{2,3\}$-group.
Furthermore, Theorem A answers a question of Babai [1, Question 7d], at least for 
finite groups.

\medskip

\noindent \textbf{Theorem A} Let $(G, \Omega)$ be a finite solvable permutation
group. There exists a subset $\Delta$ of $\Omega$ whose setwise stabilizer $S$
has the property that $O^2(S)$ is a (possibly trivial) elementary abelian 3-group.
In particular, $S$ is a monomial group, all of whose irreducible character degrees
are powers of 2.

\medskip

Set-stabilizers arise naturally in representation theory, when one examines the centralizers 
of vectors in imprimitive modules for finite groups. This happens, for example, in the
proof of Brauer's height zero conjecture for $p$-solvable groups [4,5]. The book [7] 
devotes a section to set-stabilizers, and applies that section to various situations later 
in the book.

The main tool in the proof of Theorem A is the classification of the primitive
solvable permutation groups with no regular orbit on the power set [2, Theorem 1].
We remark that the corresponding classification for nonsolvable groups was
achieved by Seress [9].

An interesting recent preprint of Sabatini [8]  shows that that set-stabilizers of
derived length at most 3 must exist when $G$ is solvable. He gives an example
to show that metabelian set-stabilizers need not exist. Somewhat analogously,
our Theorem A shows that set-stabilizers of Fitting height at most 2 must exist.
In Example 3.3, we show that nilpotent stabilizers need not exist. Sabatini [8, p. 9]
has now improved Example 3.3 by showing that even supersolvable set-stabilizers
need not exist when $G$ is solvable.

I thank Luca Sabatini for helpful communications.

\bigskip

\begin{center}
\textbf{2. Primitive Groups}
\end{center}

\bigskip

Let $(G,\Omega)$ be a finite solvable permutation group. If $\Gamma\leq
\Omega$, then by the stabilizer of $\Gamma$, we mean the setwise stabilizer,
sometimes denoted $\mathrm{Stab}_G(\Gamma)$. If $S\leq G$, we say that 
$S$ has the ``required structure'' if $O^2(S)$ is a possibly trivial elementary
abelian 3-group. We say that $(G,\Omega)$ is ``nice'' if there exist subsets
$\Delta_1$ and $\Delta_2$ of $\Omega $ with $|\Delta_1|<|\Delta_2|\leq
|\Omega/2|$ such that the stabilizers in $G$ of $\Delta_1$ and $\Delta_2$
have the required structure; here $\Delta_1$ may be the empty set.

In this section, we assume that $G$ is a primitive group on $\Omega$. We denote 
by $M$ a point stabilizer, so that $M$ is a maximal subgroup of $G$. Then $G$
has a regular normal elementary abelian subgroup $K$, which is self-centralizing
in $G$. We have $G=MK, M\cap K=1$, and $M$ acts faithfully and irreducibly on
$K$ by conjugation. The conjugation action of $M$ on $K$ is permutation isomorphic
to the action of $M\leq G$ on $\Omega$. We may consider $M$ as a subgroup of
$\mathrm{GL}_m(p)$, where $|K|=p^m$ for some prime number $p$. 

We will frequently encounter the semilinear group $\mathrm{S}(p^m)\leq
\mathrm{GL}_m(p)$, which is generated by the $p^m-1$ scalar multiplication maps
on the additive group of the finite field $\mathrm{GF}(p^m)$ and the $m$ field
automorphisms of $\mathrm{GF}(p^m)$. We have $|\mathrm{S}(p^m)|=m(p^m-1)$.
Now $\mathrm{S}(p^m)$ acts on the additive group of translations of $\mathrm{GF}(p^m)$.
We denote the semidirect product, the affine semilinear group, by $\mathrm{AS}(p^m)$.
We have $|\mathrm{AS}(p^m)|=m(p^m-1)p^m$.

\medskip

\noindent \textbf{Lemma 2.1} Let $(G,\Omega)$ be a primitive solvable permutation group. Let
$G_0$ be the set of elements of odd prime order in $G$. Let  $\mathcal{S}$ be the set
of all ordered pairs $(g, \Gamma)$, where $g \in G_0, \Gamma \leq \Omega$,
and $g$ stabilizes $\Gamma$. Then

\begin{enumerate}
\item We have $|\mathcal{S}|\leq 24^{-1/3}(2^{5|\Omega|/9})|\Omega|^{13/4}$.
\item If $|\mathcal{S}|<2^{|\Omega|}$, then there exists $\Gamma \leq \Omega$
such that $\mathrm{Stab}_G(\Gamma)$ is a 2-group.
\item If $|\mathcal{S}|<2^{|\Omega|-1}$ and $|\Omega|>9$, then $(G, \Omega)$
is nice.

\end{enumerate}
\textit{Proof.} Let $g\in G_0$. The number of subsets of $\Omega$ stabilized by $g$ is
$2^{c(g)}$, where $c(g)$ is the nuber of cycles of $g$ on $\Omega$. We claim that
$c(g) \leq 5|\Omega|/9$. To see this, let $|\Omega|=|K|=q^m$ for a prime number $q$
and a positive integer $m$. Let $g$ have order $p$. Then
$$c(g)=|C_K(g)|+|K-C_K(g)|/p.$$
If $q\ne 2$, then
$$c(g)\leq q^{m-1}+(q^m-q^{m-1})/p \leq q^{m-1}+(q^m-q^{m-1})/3.$$
Hence 
$$c(g) \leq q^m(1/q+(q-1)/3q) \leq 5q^m/9=5|\Omega|/9.$$
If $q=2$, then $c(g)\leq 2^{m-2}+(2^m-2^{m-2})/p$, which is at most
$2^m(1/4+(3/4)(1/3))=|\Omega|/2.$ This proves the claim. By Wolf's bound
([10, Corollary 3.3], [7, Corollary 3.6]), we have $|G_0|<|G|\leq 24^{-1/3}
|\Omega|^{13/4}$. Hence
$$|\mathcal{S}|\leq|G_0|2^{5|\Omega|/9}\leq (24^{-1/3})(2^{5|\Omega|/9})|\Omega|^{13/4},$$
proving (1).

If $|\mathcal{S}|<2^{|\Omega|},$ then there exists $\Gamma \leq \Omega$ such that $\mathcal S$
contains no ordered pair $(g, \Gamma)$. Hence $\mathrm{Stab}_G(\Gamma)$ is a 
2-group, proving (2). To prove (3), suppose, to get a contradiction, that there is an integer
$k$ such that every subset of $\Omega$ whose stabilizer in $G$ is a 2-group has 
cardinality $k$ or $|\Omega|-k$. Then every subset of $\Omega$ whose cardinality is not 
$k$ or $|\Omega|-k$ is stabilized by an element of $G_0$ and so $|\mathcal{S}| \geq 2^{|\Omega|}
-2{|\Omega| \choose k}$. Now an easy induction argument using Pascal's triangle shows
that if $n \geq 9$ and $0 \leq r \leq s \leq n$, then ${n \choose r}+{n \choose s}\leq 2^{n-1}$.
Thus $|\mathcal{S}| \geq 2^{|\Omega|-1},$ a contradiction. This proves (3).

\medskip

\noindent \textbf{Lemma 2.2} If $(G, \Omega)$ is primitive and $|\Omega| \geq 49$ or if
$|\Omega|=25$ or $|\Omega|=32$, then $(G, \Omega)$ is nice.

\medskip
\noindent \textit{Proof.} Suppose $|\Omega| \geq 49.$ By Lemma 2.1, it suffices to show that
$$(24^{-1/3})|\Omega|^{13/4}(2^{5|\Omega|/9})<2^{|\Omega|-1}.$$
This is equivalent to
$$(13/4)\mathrm{log}_2 |\Omega|-(1/3)\mathrm{log}_2(24)<(4/9)|\Omega|-1,$$
which holds for $|\Omega| \geq 49.$

 Suppose next that $|\Omega|=25.$ By [7, Theorem 2.11], $M$ is isomorphic to a subgroup
 of $\mathrm{GL}(2, 5)$ of order dividing 96. Elements of order 3 in $G$ fix only the zero
 vector in $K$, and so have 9 cycles on $\Omega$. By Sylow's theorem, $G$ has at most
 $32 \cdot 25$ subgroups of of order 3. Hence the number of subsets of $\Omega$ stabilized 
 by some element of order 3 is at most $32 \cdot 25 \cdot 2^9$. Elements of order 5 in $G$
 have 5 cycles on $\Omega$ and together stabilize at most $24 \cdot2^5$ subsets. Thus the
 number of subsets of $\Omega$ stabilized by some element of odd prime order is less than 
 $32 \cdot 25 \cdot 2^9+24 \cdot 2^5$, which is less than $2^{24}$. Hence $(G, \Omega)$ is 
 nice by Lemma 2.1
 
 Now suppose that $|\Omega|=32.$ By [7, p. 88], $G \leq \mathrm{AS}(32)$, the affine
 semilinear group over $\mathrm{GF}(32)$.  We may assume $G=\mathrm{AS}(32)$, which
 has order $2^5 \cdot 5 \cdot 31.$ Elements of order 5 have 2 fixed points in $\Omega$ and
 elements of order 31 have one fixed point. Hence no set of size 3 or 4 is stabilized by any
 element of $G_0,$ proving niceness.
 
 \medskip
 
 \noindent \textbf{Lemma 2.3} If $(G,\Omega)$ is primitive and $|\Omega|=16$ or
 $|\Omega|=27$, then $(G,\Omega)$ is nice.
 
 \medskip
 \noindent \textit{Proof.} Let $G=KM$ as above. We may assume that $M$ is a maximal
 irreducible solvable subgroup of $\mathrm{GL}(K)$. By [7, pp. 89-90], we have four possibilities:
 
 \begin{enumerate}
 \item $|\Omega|=16, |M|=60, G=\mathrm{AS}(16)$
 \item $|\Omega|=16, M\cong S_3 \wr Z_2, G\cong S_4 \wr Z_2$
 \item $|\Omega|=27, |M|=78, G=\mathrm{AS}(27)$
 \item $|\Omega|=27, M\cong Z_2 \wr S_3, G\cong S_3 \wr S_3$
 \end{enumerate}
 In case (2), the action of $G$ on $\Omega$ is permutation isomorphic to the product
 action of $S_4 \wr Z_2$ on $F\times F$, where $F=\{1,2,3,4\}$. This means that the
 base group of the wreath product acts componentwise and the outside $Z_2$ permutes 
 the coordinates. In case (4), we have the analogous product action of $S_3 \wr S_3$
 on $T\times T\times T$, where $T=\{1,2,3\}$.
 
 In case (1), a regular $G$-set of size 6 exists by [2, p. 62] and a regular $G$-set of
 size 7 exists by [7, p. 89]. Hence $(G, \Omega)$ is nice. In case (2), take $\Delta_1=
 \{(1,1)\}$ and $\Delta_2=\{(1,1),(2,2)\}$, viewed as subsets of $F \times F$. The stabilizer
 of $\Delta_1$ in the base group is isomorphic to $S_3 \times S_3$ and so 
 $O^2(\mathrm{Stab}_G(\Delta_1)) \cong Z_3 \times Z_3$. The stabilizer in $G$ of
 $\Delta_2$ is a 2-group. Hence $(G, \Omega)$ is nice.
 
 In case (3), there is a regular $G$-set of size 11 by [7, p. 90] and there is a regular
 $G$-set of size 4 by [2, p. 63]. Hence $(G, \Omega)$ is nice. In case (4), we note
 that $|O^2(G)|=648$. By Lemma 2.1, each element of order 3 has at most 15 cycles on $\Omega$. 
 Hence elements of order 3 stabilize at most $648 \cdot 2^{15}<2^{26}$ subsets of
 $\Omega$. Thus $(G,\Omega)$ is nice by Lemma 2.1.
 
 \medskip
 
 \noindent \textbf{Lemma 2.4} Suppose $|\Omega|$ is a prime number $p \geq 5$. Then 
 $(G, \Omega)$ is nice.
 
 \medskip
 \noindent \textit{Proof.} Here $G$ is a Frobenius group on $\Omega$, and so only the
 identity fixes two points. Hence if we choose $\Delta_1 \leq \Omega$ with $|\Delta_1|
 =2$, then $\mathrm{Stab}_G(\Delta_1)$ has order 1 or 2. Now if 3 does not divide $p-1$,
 then no element of odd prime order stabilizes any 3-subset of $\Omega$. Hence $(G, \Omega)$
 is nice when 3 does not divide $p-1$ and $p>5$. If $p=5$ take $|\Delta_1|=1$ and $|\Delta_2|=2$.
 
 Next suppose that 3 divides $p-1$. Then $G$ has exactly $p$ subgroups of order 3, each of which
 stabilizes $(p-1)/3$ three-element subsets of $\Omega$. Hence the number of 3-subsets of
 $\Omega$ stabilized by some element of order 3 in $G$ is $p(p-1)/3$, which is less than 
 $p(p-1)(p-2)/6$. Thus if we take $\Delta_2$ to be a 3-subset of $\Omega$ that is stabilized by no
 element of order 3 in $G$, then  $\mathrm{Stab}_G(\Delta_2)$ has order 1 or 2, and so
 $(G, \Omega)$ is nice.
 
 \medskip
 \noindent \textbf{Corollary 2.5} If $(G, \Omega)$ is primitive and $|\Omega|>9$, then 
 $(G, \Omega)$ is nice.
 
 \medskip
 \noindent \textit{Proof.} This follows from Lemmas 2.2, 2.3, and 2.4.
 
 \medskip
 \noindent \textbf{Lemma 2.6} Suppose that $(G, \Omega)$ is primitive and $2 \leq |\Omega|
 \leq 9$.Then $(G, \Omega)$ is nice.
 
 \medskip
 \noindent \textit{Proof.} We have $G=MK$ as above, with $K=|\Omega|$. Since $|K| \leq 9$, 
 there is a unique maximal irreducible solvable linear group on $K$. The cases $|\Omega|=$
 2, 3, and 4 are trivial; we take $(|\Delta_1|,|\Delta_2|)=$ (0,1), (0,1), and (1,2) respectively.
 The cases $|\Omega|=5$ and $|\Omega|=7$ were covered in Lemma 2.4.
 
 Suppose $|\Omega|=8$. Then $\mathrm{S}(8)$ is the unique maximal irreducible solvable
 subgroup of $\mathrm{GL}(3,2)$ by [7, Corollary 2.13]. Thus we may assume that $G=\mathrm{AS}(8)$.
 Choose $\Delta_1 \leq \Omega$ with $|\Delta_1|=2$. Since $G$ is doubly transitive on $\Omega, G$
 is transitive on the 2-subsets of $\Omega$. There are 28 such subsets and so $\mathrm{Stab}_G(\Delta_1)$
 has order 6, and must therefore have the required structure. Choose $\Delta_2 \leq \Omega$ with
 $|\Delta_2|=3.$ Since involutions in $G$ are fixed point free on $\Omega$, no involution in $G$
 stabilizes $\Delta_2.$ Obviously no element of order 7  stabilizes $\Delta_2$. Since $\Omega$ has
 56 subsets of size 3, it follows that $\mathrm{Stab}_G(\Delta_2) $ is nontrivial, and hence has order 3. Thus $(G, \Omega)$
 is nice with $(|\Delta_1|, |\Delta_2|)=(2,3)$.
 
 Now suppose that $|\Omega|=9$. Then $G$ is contained in the affine general linear group
 $\mathrm{AGL}(2,3)$. Since $\mathrm{ASL}(2,3)$ is an index 2 subgroup of $\mathrm{AGL}(2,3)$,
 the definition of niceness allows us to assume that $G=\mathrm{ASL}(2,3)$, which has order 216.
 Now $G$ is doubly transitive on $\Omega$, and therefore transitive on the 2-subsets of $\Omega$.
 Since there are 36 such subsets, the stabilizer of a fixed 2-element subset $\Delta_1$ has order 6,
 and must therefore have the required structure.
 
 Let $t \in \mathrm{SL}(2,3) \leq G$ be a fixed element of order 3. Let $\Delta_2 \leq \Omega$ consist 
 of one fixed point of $t$ together with three points that are cyclically permuted by $t$. We claim that
 $S:=\mathrm{Stab}_G(\Delta_2)$ is equal to $\langle t \rangle$. To see this, first note that elements
 of $O_3(G)$ are fixed point free on $\Omega$, and so $S \cap O_3(G)=1$ and $S$ is isomorphic to
 a subgroup of $\mathrm{SL}(2,3)$. On the other hand, no nonidentity linear transformation on $K$
 can fix more than three vectors in $K$, and so no nonidentity element of $G$ can fix more than
 three points in $\Omega$. Thus $S$ acts faithfully on $\Delta_2$, and so $S$ is also isomorphic
 to a subgroup of $S_4$. This is possible only if $S=\langle t \rangle$. Hence $(G, \Omega)$ is 
 nice with $|\Delta_1|=2$ and $|\Delta_2|=4$.
 
 \medskip
 
 \noindent \textbf{Lemma 2.7} In the situation of Lemma 2.6, suppose that $\Delta$ is a
 3-element subset of $\Omega$. If $|\Omega|=7$ or $|\Omega|=8$, then no nonidentity 
 element of $G$ fixes $\Delta$ pointwise. If $|\Omega|=9$, then there exists a 3-element 
 subset of $\Omega$ whose pointwise stabilizer in $G$ is trivial.
 
 \medskip
 \noindent \textit{Proof.} This is obvious when $|\Omega|=7$, because $(G, \Omega)$ is
 a Frobenius group. When $|\Omega|=8,$ a point stabilizer in $G$ is a Frobenius group
 on 7 points., so the pointwise stabilizer of any 3-subset of $\Omega$ is trivial. When 
 $|\Omega|=9$, we may assume that $G=\mathrm{AGL}(2,3)$. We identify $\Omega$
 with $V$, a 2-dimensional vector space over $\mathrm{GF}(3)$. Let $\{v,w\}$ be a basis
 of $V$. Let $\Delta=\{0,v,w\}$. Any $g \in G$ that fixes $\Delta$ pointwise fixes 0, and
 so $g \in \mathrm{GL}(2,3)$. Since $g$ also fixes $v$ and $w$, we have $g=1$.
 
 \bigskip
 \begin{center}
 \textbf{3. The General Case}
 \end{center}
 \medskip
 
 In this section, we no longer assume that $(G, \Omega)$ is primitive, but ``nice'' and
 ``required structure'' have the same meaning as before.
 
 \medskip
 \noindent \textit{Notation.} Let $(G_1, \Omega_1)$ be a permutation group and let $W$
 be a transitive permutation group on $\{1,\ldots,n\}$. The wreath product $G_1 \wr W$
 acts imprimitively on the disjoint union $\Omega=\Omega_1 \cup \ldots \cup \Omega_n$,
 where each $|\Omega_i|=|\Omega_1|$. Let $B \cong {G_1}^{n}$ denote the base group
 of this wreath product. Let $\{1=w_1,w_2, \ldots,w_n\}$ be a transversal to the stabilizer 
 in $W$ of $1 \in \{1,\ldots,n\}$, with $w_i(\Omega_1)=\Omega_i$.
 
 If $\Delta_1, \ldots ,\Delta_n$ are (not necessarily distinct) subsets of $\Omega_1$, we
 denote by $(\Delta_1, \ldots , \Delta_n)$ the disjoint union $\Delta_1 \cup w_2(\Delta_2)
 \ldots \cup w_n(\Delta_n)$. Of course $|w_i(\Delta_i)|=|\Delta_i|$ for $1 \leq i \leq n$.
 The stabilizer in $B$ of $(\Delta_1, \ldots, \Delta_n)$ is
 isomorphic to the direct product $\mathrm{Stab}_{G_1}(\Delta_1) \times \ldots \times \mathrm{Stab}
 _{G_1}(\Delta_n)$.
 
 \medskip
 \noindent \textbf{Theorem 3.1} Let $(G, \Omega)$ be a primitive solvable permutation group.
 If $|\Omega|>9,$ then there exists a subset $\Delta$ of $\Omega$ with $|\Delta| \ne |\Omega|/2$
 such that $\mathrm{Stab}_G(\Delta)=1$.
 
 \medskip
 \noindent \textit{Proof.} This is a less detailed statement of [2, Theorem 1], which is reproved
 as [7, Theorem 5.6].
 
 \medskip
 \noindent \textbf {Theorem 3.2} Let $(G, \Omega)$ be a solvable permutation group with 
 $G>1.$ Then $(G, \Omega)$ is nice.
 
 \medskip
 \noindent \textit{Proof.} We proceed by induction on $|\Omega|$. By an obvious semidirect
 product argument, we may assume that $G$ is transitive on $\Omega.$  By Corollary 2.5
 and Lemma 2.6, we may assume that $G$ is imprimitive on $\Omega$. Then $G \leq G_1 \wr W$ 
 as above. We may assume that $W$ acts primitively on $\{1, \ldots, n\}$ with $n>1$ and
 $|G_1|>1$. Since the property that a group $S$ has the required structure is inherited by
 subgroups of $S$, we may assume that $G=G_1 \wr W$.
 
 Suppose that $n>9.$ By Theorem 3.1, there exists $\Gamma \leq\{1, \ldots , n \}$ such that
 $|\Gamma| \ne n/2$ and $\mathrm{Stab}_W(\Gamma)=1$. After renumbering, we may assume that
 $\Gamma=\{1,2, \ldots, k \}$ with $0<k<n/2$. By the inductive hypothesis, there exist subsets
 $\Delta_1$ and $\Delta_2$ of $\Omega_1$ such that $\mathrm{Stab}_{G_1}(\Delta_j)$ has
 the required structure for $j=1,2,$ and $|\Delta_1| <|\Delta_2| \leq |\Omega_1|/2.$ Define
 $\widetilde{\Delta_1}$ to be the $n$-tuple $(\Delta_2, \ldots , \Delta_2, \Delta_1, \ldots , \Delta_1)$,
 with the first $k$ entries $\Delta_2$ and the remaining entries $\Delta_1$. Define $\widetilde
{ \Delta_2 }$ similarly, with the roles of $\Delta_1$ and $\Delta_2$ reversed. Then
$|\widetilde{\Delta_j}| \leq |\Omega|/2$ for $j=1,2$. Since $k<n/2,$ we have 
$|\widetilde{\Delta_1}|<|\widetilde{\Delta_2}|$. Since $|\Delta_1| \ne |\Delta_2|$, and $\Gamma$ 
has trivial stabilizer in $W$, it follows that $\mathrm{Stab}_G(\widetilde{\Delta_1})=
\mathrm{Stab}_B(\widetilde{\Delta_1})$, a direct product of groups with the required structure.
The same is true of $\mathrm{Stab}_G(\widetilde{\Delta_2})$. Thus $(G, \Omega)$ is nice if
$n>9$.

Now suppose that $n\leq 9$. Again by the inductive hypothesis there exist subsets $\Delta_1$
and $\Delta_2$ of $\Omega_1$, with $|\Delta_1|<|\Delta_2| \leq |\Omega_1|/2$, such that the stabilizers
in $G_1$ of $\Delta_1$ and $\Delta_2$ have the required structure. If $n=2$, let $\widetilde{\Delta_1}=
(\Delta_1, \Delta_1)$ and $\widetilde{\Delta_2}=(\Delta_2, \Delta_2)$. Since $|G:B|=2$, the conditions for 
niceness are easily checked. If $n=3$, let $\widetilde{\Delta_1}=(\Delta_1, \Delta_1, \Delta_2)$ and
let $\widetilde{\Delta_2}=(\Delta_1, \Delta_2, \Delta_2)$.Then $|\widetilde{\Delta_1}|<
|\widetilde{\Delta_2}| \leq |\Omega|/2$. No element of order 3 in $W \leq S_3$ can stabilize 
$\widetilde{\Delta_j}$ for $j=1,2$. Hence $\mathrm{Stab}_B(\widetilde{\Delta_j})$ has index at most 2
in $\mathrm{Stab}_G(\Delta_j)$ for $j=1,2.$ Niceness now follows as in the $n=2$ case. A similar
argument works when $n=5$; take $\widetilde{\Delta_1}=(\Delta_1, \Delta_1, \Delta_1, \Delta_2,\Delta_2)$
and $\widetilde{\Delta_2}=(\Delta_1, \Delta_1, \Delta_2, \Delta_2, \Delta_2)$.

If $n=4$, we may assume that $W=S_4$. Let $\widetilde{\Delta_1}=(\Delta_1, \Delta_1, \Omega_1-
\Delta_1, \Delta_2)$ and $\widetilde{\Delta_2}=(\Delta_1, \Delta_2, \Omega_1-\Delta_1, \Delta_2)$.
Since $\Delta_1, \Delta_2,$ and $\Omega_1-\Delta_1$ have distinct cardinalities, no element of
order 3 in $W$ can stabilize $\widetilde{\Delta_1}$ or $\widetilde{\Delta_2}$. Arguing as above, we see
that the stabilizers in $G$ of the last two sets have the required structure. Furthermore, $|\widetilde
{\Delta_1}|<|\widetilde{\Delta_2}|\leq |\Omega_1|+|\Omega_1|=|\Omega|/2.$ Hence$(G, \Omega)$
is nice.

The same idea works in the remaining cases. If $n=7$, let $$\widetilde{\Delta_1}=
(\Delta_2, \Delta_2, \Omega_1-\Delta_1, \Delta_1, \Delta_1, \Delta_1, \Delta_1)$$
and let $$\widetilde{\Delta_2}=(\Delta_1, \Delta_1, \Omega_1-\Delta_1, \Delta_2,
\Delta_2, \Delta_2, \Delta_2).$$ Since $\Delta_1, \Delta_2,$ and $\Omega_1-\Delta_1$
have distinct cardinalities, Lemma 2.7 implies that no element of order 3 in $W$ 
stabilizes $\widetilde{\Delta_1}$ or $\widetilde{\Delta_2}$. Hence $\mathrm{Stab}_G
(\widetilde{\Delta_j})/\mathrm{Stab}_B(\widetilde{\Delta_j})$ is a 2-group for $j=1,2.$
One checks that $|\widetilde{\Delta_1}|<|\widetilde{\Delta_2}|<|\Omega|/2.$ Hence 
$(G, \Omega)$ is nice.

For $n=8$, we may assume $W=\mathrm{AS}(8)$. We append a final $\Delta_1$ to
the $\widetilde{\Delta_1}$ of the preceding paragraph, and a final $\Delta_2$ to the
$\widetilde{\Delta_2}$ of the preceding paragraph. As above, $(G, \Omega)$ is nice.

If $n=9,$ we renumber to ensure that $\{1,2,3\}$ is fixed pointwise by no nonidentity element of 
$W=\mathrm{AGL}(2,3)$. We append a final $\Delta_1$ to the $\widetilde{\Delta_1}$ of the
preceding paragraph and a final $\Delta_2$ to the $\widetilde{\Delta_2}$ of the preceding
paragraph. As above, $(G, \Omega)$ is nice.

\medskip
\noindent \textbf{Theorem A} Let $(G, \Omega)$ be a finite solvable permutation group.
There exists a subset $\Delta$ of $\Omega$ whose setwise stabilizer $S$ has the property
that $O^2(S)$ is a possibly trivial elementary abelian 3-group. In particular, $S$ is a monomial
group, all of whose irreducible character degrees are powers of 2.

\medskip
\noindent \textit{Proof.} The first statement is immediate from Theorem 3.2. The other statements
are well-known consequences of the first; see e.g. [6, Theorem 6.15, Theorem 6.23].

\medskip
\noindent \textbf{Example 3.3}  Theorem A implies that the minimum Fitting height of a 
set-stabilizer in $(G,\Omega)$ is at most two. It is easy to see that we cannot improve this to
Fitting height one; i. e. nilpotent set-stabilizers need not exist. Indeed let $G=S_4 \wr S_4$,
acting imprimitively on a set $\Omega$ of cardinality 16. Thus $\Omega$ is the disjoint union
of four sets $\Omega_i, 1\leq i\leq4$, transitively permuted by $G$. For $\Delta \leq \Omega$, 
let $\Delta_i=\Delta \cap \Omega_i$. If some $|\Delta_i|\ne 2$, then the stabilizer in the base group of $\Delta$ involves $S_3$,
and so the stabilizer in $G$ of $\Delta$ is not nilpotent. On the other hand, if all $|\Delta_i|=2$, then a 
subgroup of $G$ isomorphic to $S_4$ stabilizes $\Delta$, so again $\mathrm{Stab}_G(\Delta)$
is not nilpotent.

\newpage
\begin{center}
\textbf{References}
 \end{center}
 
 \bigskip
 
 \noindent \textbf1. L. Babai, Asymptotic coloring of locally finite graphs and profinite permutation groups,
 J. Algebra 607(2022), 64-106.
 
 \medskip
 \noindent \textbf2. D. Gluck, Trivial set-stabilizers in finite permutation groups, Canad. J. Math 35(1983), 59-67.
 
 \medskip
 \noindent \textbf3. D. Gluck and O. Manz, Prime factors of character degrees of solvable groups, 
 Bull. London Math. Soc. 19(1987), 431-437.
 
 \medskip
 \noindent \textbf4. D. Gluck and T. R. Wolf, Defect groups and character heights in blocks of 
 solvable groups, II, J. Algebra 87(1984), 222-246.
 
 \medskip
 \noindent \textbf5. D. Gluck and T. R. Wolf, Brauer's height conjecture for $p$-solvable groups,
 Trans. Amer. Math. Soc. 282(1984), 137-152.
 
 \medskip
 \noindent \textbf6. I. M. Isaacs, Character theory of finite groups, Academic Press, New York, 1976.
 \medskip
 
 \noindent \textbf7. O. Manz and T. R. Wolf, Representations of solvable groups, Cambridge
 University Press, Cambridge, 1993.
 
 \medskip 
 \noindent \textbf8. L. Sabatini, On stabilizers in finite permutation groups, 
 arXiv: 2411.18534.
 
 \medskip
 \noindent \textbf9. A. Seress, Primitive groups with no regular orbits on the
 set of subsets, Bull. London Math. Soc. 29(1997), 697-704.
 
 \medskip
 \noindent \textbf{10}. T. R. Wolf, Solvable and nilpotent subgroups of 
 GL($n,q^m)$, Canad. J. Math. 34(1982), 1097-1111.

\end{document}